\documentclass[abstract=true,version=first]{scrartcl}

\usepackage{typearea}
\typearea{22}
\usepackage{amsmath,amssymb}
\usepackage{amsthm}
\theoremstyle{plain}
\newtheorem{theorem}{Theorem}[section]
\newtheorem{lemma}[theorem]{Lemma}
\newtheorem{corollary}[theorem]{Corollary}

\theoremstyle{definition}
\newtheorem{definition}{Definition}[section]
\newtheorem{example}{Example}[section]
\theoremstyle{remark}
\newtheorem{remark}{Remark}
\usepackage[pdftex]{graphicx}
\usepackage[colorlinks,linkcolor=blue,citecolor=blue,urlcolor=blue]{hyperref}
\usepackage{mathrsfs}
\newcommand{\mathbd}[1]{\boldsymbol{#1}}
\newcommand{\domD}{\mathscr{D}}
\newcommand{\dd}{\,\mathrm{d}}
\newcommand{\pim}{\pi_{\mathrm{m}}}

\DeclareMathOperator{\OO}{O}
\DeclareMathOperator{\sinc}{sinc}

\DeclareMathOperator{\Beta}{B}

\renewcommand{\Im}{\operatorname{Im}}
\newcommand{\trans}{\mathrm{T}}
\newcommand{\LC}{\mathbf{L}}
\newcommand{\MC}{\mathbf{M}}
\newcommand{\HC}{\mathbf{HC}}
\newcommand{\textSE}{\scriptscriptstyle\mathrm{SE}}
\newcommand{\textDE}{\scriptscriptstyle\mathrm{DE}}
\newcommand{\SEt}{\psi^{\textSE}}
\newcommand{\DEt}{\psi^{\textDE}}
\newcommand{\SEtInv}{\{\SEt\}^{-1}}
\newcommand{\DEtInv}{\{\DEt\}^{-1}}
\newcommand{\SEtDiv}{\{\SEt\}'}
\newcommand{\DEtDiv}{\{\DEt\}'}
\newcommand{\ProjSE}{\mathcal{P}_N^{\textSE}}
\newcommand{\ProjDE}{\mathcal{P}_N^{\textDE}}
\newcommand{\uSEn}{u_N^{\textSE}}
\newcommand{\vSEn}{v_N^{\textSE}}
\newcommand{\uDEn}{u_N^{\textDE}}
\newcommand{\vDEn}{v_N^{\textDE}}

\newcommand{\xSE}{t^{\textSE}}
\newcommand{\xDE}{t^{\textDE}}
\newcommand{\tSE}{\tilde{t}^{\textSE}}
\newcommand{\tDE}{\tilde{t}^{\textDE}}
\newcommand{\Fred}{\mathcal{K}}
\newcommand{\FredSEn}{\Fred^{\textSE}_N}
\newcommand{\FredDEn}{\Fred^{\textDE}_N}
\newcommand{\eSEn}{E_{n}^{\textSE}}
\newcommand{\kSEn}{K_{n}^{\textSE}}
\newcommand{\gSEn}{\mathbd{g}_n^{\textSE}}
\newcommand{\eDEn}{E_{n}^{\textDE}}
\newcommand{\kDEn}{K_{n}^{\textDE}}
\newcommand{\gDEn}{\mathbd{g}_n^{\textDE}}
\usepackage{tikz}
\numberwithin{equation}{section}

\begin{document}

\title{Sinc-collocation methods with consistent collocation points for Fredholm integral equations of the second kind}
\author{Tomoaki Okayama\footnote{Graduate School of Information Sciences,
Hiroshima City University}}



\maketitle


\begin{abstract}
Sinc-collocation methods are known to be efficient
for Fredholm integral equations of the second kind,
even if functions in the equations have endpoint singularity.
However, existing methods have the disadvantage of
inconsistent collocation points.
This inconsistency complicates the implementation of such methods,
particularly for large-scale problems.
To overcome this drawback,
this study proposes another Sinc-collocation methods
with consistent collocation points.
The results of a theoretical error analysis show that the proposed methods
have the same convergence property as existing methods.
Numerical experiments suggest the superiority of the proposed methods
in terms of implementation and computational cost.
\end{abstract}


\section{Introduction}

This paper is concerned with
Fredholm integral equations of the second kind of the following form:
\begin{equation}
 u(t) - \int_a^b k(t, s)u(s) \dd s = g(t),
\quad a\leq t\leq b,
\label{eq:Fredholm}
\end{equation}
where $k$ and $g$ are given continuous functions, and $u$ is the solution
to be determined.
Most numerical methods provided in the literature
do not perform well when the functions $k$ and $g$
have derivative singularity at the endpoints.
To overcome the difficulty,
Rashidinia and Zarebnia~\cite{RZ1} proposed a Sinc-collocation method,
which was derived without assuming differentiability at the endpoints.
The results of numerical experiments indicated that their method
can achieve $\OO(\exp(-c\sqrt{N}))$,
where $N$ denotes approximately half the number of collocation points.
However, they did not prove whether their method was guaranteed to converge.
Furthermore, their method required information on the boundary values of
the solution $u(t)$, which are not available in practice.

To address this situation,
Okayama et al.~\cite{OMS} reformed their method
so that it could be implemented without the information of the boundary
values of $u(t)$, and proved that
the reformed method can attain $\OO(\exp(-c\sqrt{N}))$.
The method is based on the Sinc approximation combined with
the Single-Exponential (SE) transformation:
\begin{equation}
t = \SEt(x)
  = \frac{b-a}{2}\tanh\left(\frac{x}{2}\right)+\frac{b+a}{2}.
\label{eq:SE-trans}
\end{equation}
They further improved the method by replacing the SE transformation
with the Double-Exponential (DE) transformation:
\begin{equation}
 t = \DEt(x)
  = \frac{b-a}{2}\tanh\left(\frac{\pi}{2}\sinh x\right)+\frac{b+a}{2},
\label{eq:DE-trans}
\end{equation}
and proved that their improved method can attain
a higher convergence rate: $\OO(\exp(-c N/\log N))$.

The methods were derived as follows.
Based on the Sinc approximation with auxiliary functions,
set the approximate solution $u_N(t)$ as
\[
 u_N(t) = c_{-N-1} \frac{b-t}{b-a}
+ \sum_{j=-N}^{N} c_{j}
\sinc\left(\frac{\psi^{-1}(t) - jh}{h}\right)
+ c_{N+1} \frac{t-a}{b-a},
\]
where $\sinc x = \sin(\pi x)/(\pi x)$,
and $\psi$ is either $\SEt$ or $\DEt$.
Substituting $u_N$ into the given equation~\eqref{eq:Fredholm}
and setting $(2N+3)$ collocation points, $t=t_i$ $(i=-N-1,\,\ldots,\,N+1)$, as
\[
  t_i =
 \begin{cases}
  a               & (i = -N-1),\\
  \psi(ih) & (i = -N, \,\ldots,\,N),\\
  b               & (i =  N+1),
 \end{cases}
\]
we obtain a system of linear equations with respect to
unknown coefficients $c_j$ $(j=-N-1,\,\ldots,\,N+1)$.
Those collocation points are chosen
so that $u_N(t_i) = c_i$ holds for $i = -N-1,\,\ldots,\,N+1$.
However, because of the inconsistency of the collocation points at $i=\pm (N+1)$,
the methods are rather painful to implement,
particularly for \emph{systems of} Fredholm integral equations.

As a related work, Sinc-Nystr\"{o}m methods have been
proposed~\cite{MM}\footnote{Although the title of the paper~\cite{MM} refers to
Sinc-collocation methods, the work considers Sinc-Nystr\"{o}m methods
in reality.},
which were derived as follows.
Set the approximate solution $v_N(t)$ as
\[
 v_N(t)
= g(t) + h \sum_{j=-N}^N k(t, \psi(jh)) v_N(\psi(jh)) \psi'(jh),
\]
which is obtained by discretizing the integral in~\eqref{eq:Fredholm}.
Setting $(2N+1)$ collocation points $t=\tilde{t}_i$ $(i=-N,\,\ldots,\,N)$ as
\[
 \tilde{t}_i = \psi(ih) \quad (i = -N,\,\ldots,\,N),
\]
we obtain a system of linear equations with respect to
unknown coefficients $v_N(\psi(jh))$ $(j=-N,\,\ldots,\,N)$.
The collocation points $\tilde{t}_i$ are consistent
for all $i=-N,\,\ldots,\,N$, which is a considerable advantage in terms of practical implementation.
In fact, the consistency was effectively utilized
in the case of a \emph{system of} Fredholm integral equations~\cite{RZ2}.
However, the approximate solution $v_N(t)$ includes
the given functions, $g(t)$ and $k(t, s)$,
which may require a high computational cost for evaluation.

To address this issue,
this study derives another Sinc-collocation methods
based on the concept of the Sinc-collocation methods
derived for initial value problems~\cite{TO}.
The approximate solution $\tilde{u}_{N}(t)$ is given
by using $\tilde{c}_j = v_N(\psi(jh))$ as
\begin{align*}
 \tilde{u}_N(t)
&=\tilde{c}_{-N}\frac{b-t}{b-a}
+\sum_{j=-N}^{N}\left(\tilde{c}_j
 - \tilde{c}_{-N}\frac{b-\psi(jh)}{b-a}
 - \tilde{c}_{N} \frac{\psi(jh)-a}{b-a}
\right)\sinc\left(\frac{\psi^{-1}(t) - jh}{h}\right)
+\tilde{c}_{N}\frac{t-a}{b-a},
\end{align*}
where the basis functions of $\tilde{u}_N(t)$ include
only elementary functions. In addition,
the coefficients $\tilde{c}_j$ $(j=-N,\,\ldots,\,N)$
are obtained by solving the system of linear equations
of the Sinc-Nystr\"{o}m methods,
where the consistent collocation points $t=\tilde{t}_i$ are employed.
Therefore, the derived methods have an advantage in
terms of both implementation and computational cost.

Furthermore, by conducting a theoretical convergence analysis,
this study proves that the new method
with $t = \SEt(x)$
can attain $\OO(\exp(-c\sqrt{N}))$,
and the new method
with $t = \DEt(x)$
can attain $\OO(\exp(-c N/\log N))$.
These convergence rates are the same as those of the original
Sinc-collocation methods proposed by Okayama et al.~\cite{OMS}.
However, the proposed methods may require less computation time
compared to the original methods to obtain the same precision.
This is due to the difference in the system of linear equations.
This point is discussed in detail in Section~\ref{sec:discussion}.

The remainder of this paper is organized as follows.
In Section~\ref{sec:fundamental}, as a preliminary,
basic definitions and theorems for Sinc methods are explained.
Sinc methods include two important parameters, $\alpha$ and $d$,
which are used to set the mesh size, $h$.
Section~\ref{sec:smooth-sol} describes how to find those parameters,
based on the result on the smoothness property
of the solution of the equation~\eqref{eq:Fredholm}.
In Section~\ref{sec:original-Sinc},
the original Sinc-collocation methods are described.
In Section~\ref{sec:new-Sinc},
the newly proposed Sinc-collocation methods are described,
which are the main result of this paper.
In Section~\ref{sec:numer}, some numerical examples are presented.
In Section~\ref{sec:proofs},
convergence theorems of the new Sinc-collocation methods
are proved.
Section~\ref{sec:conclusion} presents a conclusion of this paper.

\section{Fundamental definitions and theorems for Sinc methods}
\label{sec:fundamental}

The Sinc approximation is a function approximation formula expressed as
\begin{equation}
 F(x)\approx \sum_{j=-N}^N F(jh)S(j,h)(x),\quad x\in\mathbb{R},
\label{eq:Sinc-approx}
\end{equation}
where $S(j,h)(x)$ denotes the so-called ``Sinc function'' defined by
\[
 S(j,h)(x) = \sinc\left(\frac{x - jh}{h}\right).
\]
Numerical methods based on the Sinc approximation are referred to as
``Sinc methods.''
For example, a quadrature formula obtained by integrating
both sides of~\eqref{eq:Sinc-approx} as
\begin{equation}
 \int_{-\infty}^{\infty}F(x)\dd x
\approx\sum_{j=-N}^N F(jh)\int_{-\infty}^{\infty}S(j,h)(x)\dd x
=h\sum_{j=-N}^N F(jh)
\label{eq:Sinc-quad}
\end{equation}
is referred to as the ``Sinc quadrature.''
Comprehensive summaries of existing Sinc methods
have been provided in several books and review papers~\cite{Ste1,Ste3,Ste2,SugiMatsu}.
This section provides
fundamental definitions and theorems for such Sinc methods as preliminaries for the remainder of this work.

\subsection{SE-Sinc approximation and SE-Sinc quadrature on the finite interval}

To apply the Sinc approximation~\eqref{eq:Sinc-approx}
to a function $f(t)$ on the finite interval $(a, b)$,
the SE transformation~\eqref{eq:SE-trans} is frequently utilized as
\begin{align}
 f(\SEt(x))&\approx \sum_{j=-N}^N f(\SEt(jh))S(j,h)(x),\quad x\in\mathbb{R},
\nonumber
\intertext{which is equivalent to}
 f(t)&\approx \sum_{j=-N}^N f(\SEt(jh))S(j,h)(\SEtInv(t)),\quad t\in (a,b).
\label{eq:SE-Sinc-approx}
\end{align}
We refer to this approximation as the ``SE-Sinc approximation.''
Similarly, to apply the Sinc quadrature~\eqref{eq:Sinc-quad}
to the integral on the finite interval $(a, b)$,
the SE transformation~\eqref{eq:SE-trans} is frequently utilized as
\begin{equation}
 \int_a^b f(t)\dd t
=\int_{-\infty}^{\infty} f(\SEt(x))\SEtDiv(x)\dd x
\approx h \sum_{j=-N}^N f(\SEt(jh))\SEtDiv(jh).
\label{eq:SE-Sinc-quad}
\end{equation}
We refer to this approximation as the ``SE-Sinc quadrature.''
The function space defined below is important in formalizing an error analysis
of these approximations.

\begin{definition}
\label{Def:LC}
Let $\alpha$ be a positive constant,
and let $\domD$ be a bounded and simply-connected domain
(or Riemann surface)
that satisfies $(a,\,b)\subset \domD$.
Then, $\LC_{\alpha}(\domD)$ denotes the family of all functions $f$
satisfying the following conditions:
(i) $f$ is analytic in $\domD$;
(ii) there exists a constant $C$ such that
for all $z$ in $\domD$,
\begin{equation}
|f(z)|\leq C |Q(z)|^{\alpha}, \label{Leq:LC-bounded-by-Q}
\end{equation}
where the function $Q$ is defined by $Q(z)=(z-a)(b-z)$.
\end{definition}
When the SE transformation is utilized,
the domain $\domD$ in Definition~\ref{Def:LC} should be
an eye-shaped region:
$\SEt(\domD_d) = \{z=\SEt(\zeta):\zeta\in\domD_d\}$
(see also an example figure \cite[Fig.~1]{OMS}),
where $\domD_d$ is a strip domain defined by
$\domD_d=\{\zeta\in\mathbb{C}:|\Im\zeta|<d\}$
for a positive constant $d$.
Then, error analyses of the above approximations
are described as follows.

\begin{theorem}[Stenger~{\cite[Theorem~4.2.5]{Ste1}}]
\label{Thm:SE-Sinc-Base}
Let $f\in\LC_{\alpha}(\SEt(\domD_d))$ for $d$ with $0<d<\pi$,
let $N$ be a positive integer, and let $h$ be selected by the formula
\begin{equation}
h=\sqrt{\frac{\pi d}{\alpha N}}. \label{Def-h-SE}
\end{equation}
Then, there exists a constant $C$ that is independent of $N$,
such that
\begin{equation*}
\max_{a\leq t\leq b}
\left|f(t)-\sum_{j=-N}^N f(\SEt(jh))S(j,h)(\SEtInv(t))\right|
\leq C \sqrt{N}\exp\left(-\sqrt{\pi d \alpha N}\right).
\end{equation*}
\end{theorem}

\begin{corollary}[Okayama et al.~{\cite[Corollary~2.8]{OMS2}}]
\label{Cor:SE-Sinc-quad}
Let $(f Q)\in\LC_{\alpha}(\SEt(\domD_d))$ for $d$ with $0<d<\pi$,
let $N$ be a positive integer,
and let $h$ be selected by the formula~\eqref{Def-h-SE}.
Then, there exists a constant $C$ that is independent of $N$,
such that
\[
\left|
  \int_a^b f(t)\dd t
 - h\sum_{j=-N}^N f(\SEt(jh))\SEtDiv(jh)
\right|
\leq C \exp\left(-\sqrt{\pi d \alpha N}\right).
\]
\end{corollary}

\begin{remark}
\label{rem:choice-h}
In the case of the Sinc quadrature,
the optimal choice of $h$ is not~\eqref{Def-h-SE},
but $h=\sqrt{2\pi d/(\alpha N)}$.
However, to implement the derived method easily,
the mesh size $h$ is chosen as~\eqref{Def-h-SE}
for both the SE-Sinc approximation and quadrature.
The same applies to the DE-Sinc approximation
and quadrature (described next).
\end{remark}

\subsection{DE-Sinc approximation and DE-Sinc quadrature on the finite interval}

The key role of the SE transformation~\eqref{eq:SE-trans} is
to map the whole real axis $\mathbb{R}$ onto the finite interval $(a, b)$.
Because the DE transformation~\eqref{eq:DE-trans} also plays the same role,
the SE transformation may be replaced by the DE transformation
in the above-mentioned approximations.
It has recently been known that such replacement enhances the performance
of Sinc methods~\cite{MS,SugiMatsu}.

For example,
by the replacement in the SE-Sinc approximation~\eqref{eq:SE-Sinc-approx},
we obtain
\begin{align*}
 f(t)&\approx \sum_{j=-N}^N f(\DEt(jh))S(j,h)(\DEtInv(t)),\quad t\in (a,b).
\end{align*}
We refer to this approximation as the ``DE-Sinc approximation.''
By the same replacement in the SE-Sinc quadrature~\eqref{eq:SE-Sinc-quad},
we obtain
\begin{equation*}
 \int_a^b f(t)\dd t
=\int_{-\infty}^{\infty} f(\DEt(x))\DEtDiv(x)\dd x
\approx h \sum_{j=-N}^N f(\DEt(jh))\DEtDiv(jh).
\end{equation*}
We refer to this approximation as the ``DE-Sinc quadrature.''
In these cases,
the domain $\domD$ in Definition~\ref{Def:LC}
should be $\DEt(\domD_d)=\{z=\DEt(\zeta):\zeta\in\domD_d\}$
(see also an example figure \cite[Fig.~2]{OMS}),
which is a Riemann surface infinitely rotating around the endpoints
\footnote{More detailed explanation/figure can be found in Tanaka et al.~\cite{tanaka09:_class_of}}.
Error analyses of the above approximations
are described as follows.

\begin{theorem}[Tanaka et al.~{\cite[Theorem~3.1]{tanaka09:_class_of}}]
\label{Thm:DE-Sinc-Base}
Let $f\in\LC_{\alpha}(\DEt(\domD_d))$ for $d$ with $0<d<\pi/2$,
let $N$ be a positive integer, and let $h$ be selected by the formula
\begin{equation}
h=\frac{\log(2 d N /\alpha)}{N}. \label{Def-h-DE}
\end{equation}
Then, there exists a constant $C$ that is independent of $N$,
such that
\begin{equation*}
\max_{a\leq t\leq b}
\left|f(t)-\sum_{j=-N}^N f(\DEt(jh))S(j,h)(\DEtInv(t))\right|
\leq C \exp\left\{\frac{-\pi d N}{\log(2 d N/\alpha)}\right\}.
\end{equation*}
\end{theorem}

\begin{corollary}[Okayama et al.~{\cite[Corollary~2.9]{OMS2}}]
\label{Cor:DE-Sinc-quad}
Let $(fQ)\in\LC_{\alpha}(\DEt(\domD_d))$ for $d$ with $0<d<\pi/2$,
let $N$ be a positive integer,
and let $h$ be selected by the formula~\eqref{Def-h-DE}.
Then, there exists a constant $C$ that is independent of $N$,
such that
\[
\left|\int_a^b f(t)\dd t - h\sum_{j=-N}^N f(\DEt(jh))\DEtDiv(jh)\right|
\leq C \exp\left\{\frac{-2\pi d N}{\log(2 d N/\alpha)}\right\}.
\]
\end{corollary}

After the replacement, the convergence rate is improved to
$\OO(\exp(-c N/\log N))$,
which is substantially higher than the previous rate of
$\OO(\exp(-c \sqrt{N}))$.

\begin{remark}
The convergence rate $\OO(\exp(-c N/\log N))$
is not always attainable under the assumption
in Theorem~\ref{Thm:SE-Sinc-Base} or Corollary~\ref{Cor:SE-Sinc-quad}.
We should note that $\DEt(\domD_d)$ is a wider and more complicated
domain than $\SEt(\domD_d)$; that is,
the condition imposed on $f$ in the case of the DE transformation
is more restricted than that in the case of the SE transformation.
If the assumption of Theorem~\ref{Thm:DE-Sinc-Base}
or Corollary~\ref{Cor:DE-Sinc-quad} is not fulfilled,
the SE transformation may perform better than the DE transformation.
\end{remark}

\subsection{Generalized SE/DE-Sinc approximation on the finite interval}

According to Theorems~\ref{Thm:SE-Sinc-Base}
and~\ref{Thm:DE-Sinc-Base},
the function $f$ should belong to $\LC_{\alpha}(\domD)$.
This may seem impractical because
$f(t)$ must be zero at the endpoints, $t=a$ and $t=b$, owing to
the inequality~\eqref{Leq:LC-bounded-by-Q}.
In reality, it suffices to consider the following function spaces.

\begin{definition}
Let $\domD$ be a bounded and simply-connected domain
(or Riemann surface).
Then, $\HC(\domD)$ denotes the family of all functions
that are analytic in $\domD$ and continuous on $\overline{\domD}$.
\end{definition}

\begin{definition} \label{Def:MC}
Let $\alpha$ be a constant with $0<\alpha\leq 1$
and let $\domD$ be a bounded and simply-connected domain
(or Riemann surface)
that satisfies $(a,\,b)\subset \domD$.
Then, the space $\MC_{\alpha}(\domD)$
consists of all functions $f$
satisfying the following conditions:
(i) $f\in\HC(\domD)$;
(ii) there exists a constant $C$
for all $z$ in $\domD$ such that
\begin{align*}
|f(z)-f(a)|&\leq C |z-a|^{\alpha},\\
|f(b)-f(z)|&\leq C |b-z|^{\alpha}.
\end{align*}
\end{definition}

For functions $f\in\MC_{\alpha}(\domD)$,
using auxiliary functions
\[
 \omega_a(t)=\frac{b-t}{b-a},\quad
 \omega_b(t)=\frac{t-a}{b-a},
\]
the following approximations have been proposed:
\begin{align}
f(t)&\approx\ProjSE[f](t)\nonumber\\
&=f(\SEt(-Nh))\omega_a(t) + f(\SEt(Nh))\omega_b(t) \nonumber\\
&\quad +\sum_{j=-N}^N
\left\{f(\SEt(jh)) - f(\SEt(-Nh))\omega_a(\SEt(jh))
 - f(\SEt(Nh))\omega_b(\SEt(jh))\right\}S(j,h)(\SEtInv(t)),
\label{eq:general-SE-Sinc}\\
f(t)&\approx\ProjDE[f](t)\nonumber\\
&= f(\DEt(-Nh))\omega_a(t) + f(\DEt(Nh))\omega_b(t) \nonumber\\
&\quad +\sum_{j=-N}^N
\left\{f(\DEt(jh)) - f(\DEt(-Nh))\omega_a(\DEt(jh))
 - f(\DEt(Nh))\omega_b(\DEt(jh))\right\}S(j,h)(\DEtInv(t)).
\label{eq:general-DE-Sinc}
\end{align}
We refer to~\eqref{eq:general-SE-Sinc}
as the ``generalized SE-Sinc approximation,''
and~\eqref{eq:general-DE-Sinc}
as the ``generalized DE-Sinc approximation.''
Error analyses of these approximations
are described as follows.
Here, $\|\,\cdot\,\|_{C([a,b])}$ denotes the usual uniform norm
over $[a, b]$.

\begin{theorem}[Okayama~{\cite[Theorem 3]{TO2}}]
\label{Thm:general-SE-Sinc}
Let $f\in\MC_{\alpha}(\SEt(\domD_d))$ for $d$ with $0<d<\pi$,
let $N$ be a positive integer,
and let $h$ be selected by the formula~\eqref{Def-h-SE}.
Then, there exists a constant $C$ that is independent of $N$,
such that
\[
\|f - \ProjSE f\|_{C([a,\,b])}
\leq C \sqrt{N}\exp\left(-\sqrt{\pi d \alpha N}\right).
\]
\end{theorem}

\begin{theorem}[Okayama~{\cite[Theorem 6]{TO2}}]
\label{Thm:general-DE-Sinc}
Let $f\in\MC_{\alpha}(\DEt(\domD_d))$ for $d$ with $0<d<\pi/2$,
let $N$ be a positive integer,
and let $h$ be selected by the formula~\eqref{Def-h-DE}.
Then, there exists a constant $C$ that is independent of $N$,
such that
\[
\|f - \ProjDE f\|_{C([a,\,b])}
\leq C \exp\left\{\frac{-\pi d N}{\log(2 d N/\alpha)}\right\}.
\]
\end{theorem}

\begin{remark}
Although Okayama et al.\ also used the symbols
$\ProjSE$ and $\ProjDE$ in their paper~\cite[(2.29) and (2.31)]{OMS},
their definitions differ from those in the present
paper:~\eqref{eq:general-SE-Sinc} and~\eqref{eq:general-DE-Sinc}.
This is one of key differences between Okayama et al.\ and this study.
\end{remark}

\section{Smoothness property of the solution}
\label{sec:smooth-sol}

To approximate the solution $u$ according to Theorem~\ref{Thm:general-SE-Sinc}
or~\ref{Thm:general-DE-Sinc},
$u$ must belong to $\MC_{\alpha}(\SEt(\domD_d))$
or $\MC_{\alpha}(\DEt(\domD_d))$.
It should be noted that the parameters $\alpha$ and $d$,
which indicate smoothness of the function, are
used to select the mesh size $h$ in~\eqref{Def-h-SE} or~\eqref{Def-h-DE}.
In practice, however, investigating $u$ directly is not possible
because $u$ is an \emph{unknown} function to be determined.
To improve the situation,
Okayama et al.~\cite{OMS}
provided a sufficient condition for $u\in\MC_{\alpha}(\domD)$
using the known functions $g$ and $k$, as described below.

Suppose $k(z,\cdot)\in\HC(\domD)$
and $k(\cdot,w)\in\HC(\domD)$ for all $z,\,w\in\overline{\domD}$,
and let us introduce the integral operator
$\Fred:\HC(\domD)\to\HC(\domD)$
defined by
\[
\Fred [f](t) = \int_a^b k(t,s)f(s)\dd s.
\]
Note that the equation~\eqref{eq:Fredholm}
can be rewritten as $(I-\Fred)u = g$.
Then, the following theorem holds for both $\domD=\SEt(\domD_d)$
and $\domD=\DEt(\domD_d)$.

\begin{theorem}[Okayama et al.~{\cite[Theorem~4.3]{OMS}}]
\label{Thm:Sol-Exist-Unique}
Let $k(z,\cdot)\in\HC(\domD)$
and $k(\cdot,w)\in\MC_{\alpha}(\domD)$
for all $z,\,w\in\overline{\domD}$,
and let also $g\in\MC_{\alpha}(\domD)$.
Furthermore, assume that
the homogeneous equation $(I - \Fred)f = 0$
has only the trivial solution $f\equiv 0$.
Then, the equation~\eqref{eq:Fredholm}
has a unique solution $u\in\MC_{\alpha}(\domD)$.
\end{theorem}

This theorem suggests that the smoothness parameters $\alpha$ and $d$
of the solution $u$ can be found by investigating known functions $g$
and $k$.

\section{Original Sinc-collocation methods}
\label{sec:original-Sinc}

In this section,
Sinc-collocation methods derived by Okayama et al.~\cite{OMS}
are described.
These methods were derived
in accordance with the standard collocation procedure.

\subsection{Original SE-Sinc-collocation method}

First, suppose that
the assumptions in Theorem~\ref{Thm:Sol-Exist-Unique}
are fulfilled with $\domD=\SEt(\domD_d)$.
Then, set the approximate solution $\uSEn$ as
\begin{equation}
 \uSEn(t)
= u_{-N-1}\omega_a(t) + \sum_{j=-N}^N u_j S(j,h)(\SEtInv(t))
+ u_{N+1}\omega_b(t).
\label{Def:uSEn}
\end{equation}
Let us substitute $\uSEn$ into the equation~\eqref{eq:Fredholm},
and approximate the integral operator $\Fred$
by $\FredSEn$:
\begin{equation}
\FredSEn[f](t) = h \sum_{j=-N}^N k(t,\SEt(jh))f(\SEt(jh))\SEtDiv(jh).
\label{Def:FredSEn}
\end{equation}
In both~\eqref{Def:uSEn} and~\eqref{Def:FredSEn},
the mesh size $h$ is selected by the formula~\eqref{Def-h-SE}.
Finally, by setting the collocation points as
\begin{equation}
\xSE_i =
 \begin{cases}
  a        & (i=-N-1), \\
  \SEt(ih) & (i=-N,\,\ldots,\,N), \\
  b        & (i=N+1),
 \end{cases}
\label{Def:SE-Sinc-Sampling}
\end{equation}
we obtain a $(2N+3)\times(2N+3)$ system of linear equations:
\begin{align}
&  \{\omega_a(\xSE_i)-\FredSEn[\omega_a](\xSE_i)\}u_{-N-1}
+ \sum_{j=-N}^N
  \left\{\delta_{ij}- h k(\xSE_i,\xSE_j)\SEtDiv(jh)\right\}u_j
+ \{\omega_b(\xSE_i)-\FredSEn[\omega_b](\xSE_i)\}u_{N+1}\nonumber\\
&= g(\xSE_i),\quad i=-N-1,\,-N,\,\ldots,\,N,\,N+1.
\label{Def:SE-Sinc-linear-eq}
\end{align}
Let us define $n$ by $n=2N+3$,
and let $\eSEn$ and $\kSEn$ be $n\times n$ matrices defined by
\begin{align*}
  \eSEn
&= \left[
   \begin{array}{@{\,}c|ccc|c@{\,}}
   1          & 0      & \cdots &0      & 0 \\
   \hline
   \omega_a(\xSE_{-N}) & 1      &        &\OO & \omega_b(\xSE_{-N})\\
   \vdots     &        & \ddots &       & \vdots \\
   \omega_a(\xSE_N) & \OO &        &1      & \omega_b(\xSE_N) \\
   \hline
   0          & 0      & \cdots &0      & 1
   \end{array}
   \right], \\
  \kSEn
&= \left[
   \begin{array}{@{\,}l|clc|l@{\,}}
    \FredSEn[\omega_a](a)
   &\cdots
   & h k(a,\xSE_j)\SEtDiv(jh)
   &\cdots
   &\FredSEn[\omega_b](a) \\
   \hline
    \FredSEn[\omega_a](\xSE_{-N})
   &\cdots
   & h k(\xSE_{-N},\xSE_j)\SEtDiv(jh)
   &\cdots
   &\FredSEn[\omega_b](\xSE_{-N}) \\
    \multicolumn{1}{c|}{\vdots} & &
    \multicolumn{1}{c}{\vdots} & & \multicolumn{1}{c}{\vdots}\\
    \FredSEn[\omega_a](\xSE_{N})
   &\cdots
   & h k(\xSE_{N},\xSE_j)\SEtDiv(jh)
   &\cdots
   &\FredSEn[\omega_b](\xSE_N) \\
   \hline
    \FredSEn[\omega_a](b)
   &\cdots
   & h k(b,\xSE_j)\SEtDiv(jh)
   &\cdots
   &\FredSEn[\omega_b](b)
   \end{array}
   \right].
\end{align*}
Furthermore, let
$\gSEn=[g(a),\,g(\xSE_{-N}),\,\ldots,\,g(\xSE_N),\,g(b)]^{\trans}$.
Then, the resulting system of linear equations~\eqref{Def:SE-Sinc-linear-eq}
can be written in the matrix-vector form:
\begin{equation}
(\eSEn - \kSEn)\mathbd{u}_n = \gSEn,
\label{Def:SE-Sinc-linear-eq-matvec}
\end{equation}
where $\mathbd{u}_n=[u_{-N-1},\,u_{-N},\,\ldots,\,u_N,\,u_{N+1}]^{\trans}$.
By solving this for the coefficients $\mathbd{u}_n$,
the approximate solution $\uSEn$ is determined by~\eqref{Def:uSEn}.
This is the original SE-Sinc-collocation method.
Its convergence was analyzed as follows.

\begin{theorem}[Okayama et al.~{\cite[Theorem~6.4]{OMS}}]
\label{Thm:converge_RandZ}
Suppose that
the assumptions in Theorem~\ref{Thm:Sol-Exist-Unique} are fulfilled
with $\domD=\SEt(\domD_d)$ for $d\in(0,\pi)$.
Then, there exists a positive integer $N_0$ such that
the equation~\eqref{Def:SE-Sinc-linear-eq-matvec}
is uniquely solvable for all $N\geq N_0$.
Furthermore, there exists a constant $C$ independent of $N$
such that for all $N\geq N_0$,
\[
\|u-\uSEn\|_{C([a,b])} \leq C \sqrt{N}\exp\left(-\sqrt{\pi d \alpha N}\right).
\]
\end{theorem}

\subsection{Original DE-Sinc-collocation method}

Here,
the original DE-Sinc-collocation method is described;
it replaces the SE transformation in the original SE-Sinc-collocation method
with the DE transformation.
Suppose that
the assumptions in Theorem~\ref{Thm:Sol-Exist-Unique}
are fulfilled with $\domD=\DEt(\domD_d)$.
Then,
set the approximate solution $\uDEn$ as
\begin{equation}
\uDEn(t)
= u_{-N-1}\omega_a(t) + \sum_{j=-N}^N u_j S(j,h)(\DEtInv(t)) + u_{N+1}\omega_b(t).
\label{Def:uDEn}
\end{equation}
Let us substitute $\uDEn$ into the equation~\eqref{eq:Fredholm},
and approximate the integral operator $\Fred$
by $\FredDEn$:
\begin{equation}
\FredDEn[f](t) =  h \sum_{j=-N}^N k(t,\DEt(jh))f(\DEt(jh))\DEtDiv(jh).
\label{Def:FredDEn}
\end{equation}
In both~\eqref{Def:uDEn} and~\eqref{Def:FredDEn},
the mesh size $h$ is selected by the formula~\eqref{Def-h-DE}.
Finally, by setting the collocation points as
$t=\xDE_i$ defined by
\begin{equation}
\xDE_i =
 \begin{cases}
  a        & (i=-N-1), \\
  \DEt(ih) & (i=-N,\,\ldots,\,N), \\
  b        & (i=N+1),
 \end{cases}
\label{Def:DE-Sinc-Sampling}
\end{equation}
we obtain an $n\times n$ system of linear equations (recall $n=2N+3$):
\begin{equation}
(\eDEn - \kDEn)\mathbd{u}_n = \gDEn,
\label{Def:DE-Sinc-linear-eq-matvec}
\end{equation}
where
$\gDEn=[g(a),\,g(\xDE_{-N}),\,\ldots,\,g(\xDE_N),\,g(b)]^{\trans}$,
and
$\eDEn$ and $\kDEn$ are $n\times n$ matrices defined by
\begin{align*}
  \eDEn
&= \left[
   \begin{array}{@{\,}c|ccc|c@{\,}}
   1          & 0      & \cdots &0      & 0 \\
   \hline
   \omega_a(\xDE_{-N}) & 1      &        &\OO & \omega_b(\xDE_{-N})\\
   \vdots     &        & \ddots &       & \vdots \\
   \omega_a(\xDE_N) & \OO &        &1      & \omega_b(\xDE_N) \\
   \hline
   0          & 0      & \cdots &0      & 1
   \end{array}
   \right], \\
  \kDEn
&= \left[
   \begin{array}{@{\,}l|clc|l@{\,}}
    \FredDEn[\omega_a](a)
   &\cdots
   & h k(a,\xDE_j)\DEtDiv(jh)
   &\cdots
   &\FredDEn[\omega_b](a) \\
   \hline
    \FredDEn[\omega_a](\xDE_{-N})
   &\cdots
   & h k(\xDE_{-N},\xDE_j)\DEtDiv(jh)
   &\cdots
   &\FredDEn[\omega_b](\xDE_{-N}) \\
    \multicolumn{1}{c|}{\vdots} & &
    \multicolumn{1}{c}{\vdots} & & \multicolumn{1}{c}{\vdots}\\
    \FredDEn[\omega_a](\xDE_{N})
   &\cdots
   & h k(\xDE_{N},\xDE_j)\DEtDiv(jh)
   &\cdots
   &\FredDEn[\omega_b](\xDE_N) \\
   \hline
    \FredDEn[\omega_a](b)
   &\cdots
   & h k(b,\xDE_j)\DEtDiv(jh)
   &\cdots
   &\FredDEn[\omega_b](b)
   \end{array}
   \right].
\end{align*}
By solving this for the coefficients $\mathbd{u}_n$
in~\eqref{Def:DE-Sinc-linear-eq-matvec},
the approximate solution $\uDEn$ is determined by~\eqref{Def:uDEn}.
This is the original DE-Sinc-collocation method.
Its convergence was analyzed as follows.

\begin{theorem}[Okayama et al.~{\cite[Theorem~8.2]{OMS}}]
\label{Thm:converge_DE-Sinc}
Suppose that
the assumptions in Theorem~\ref{Thm:Sol-Exist-Unique} are fulfilled
with $\domD=\DEt(\domD_d)$ for $d\in(0,\pi/2)$.
Then, there exists a positive integer $N_0$ such that
the equation~\eqref{Def:DE-Sinc-linear-eq-matvec}
is uniquely solvable for all $N\geq N_0$.
Furthermore, there exists a constant $C$ independent of $N$
such that for all $N\geq N_0$,
\[
\|u-\uDEn\|_{C([a,b])} \leq C \exp\left\{\frac{-\pi d N}{\log(2 \pi d /\alpha)}\right\}.
\]
\end{theorem}

\section{New Sinc-collocation methods with consistent collocation points}
\label{sec:new-Sinc}

In this section,
new Sinc-collocation methods are presented.
These methods are derived
based on the concept of Sinc-collocation methods
derived for initial value problems~\cite{TO},
which are closely related to the Sinc-Nystr\"{o}m methods~\cite{MM,RZ2}.
Convergence theorems for the presented methods are stated,
and proofs are provided in Section~\ref{sec:proofs}.

\subsection{New SE-Sinc-collocation method with consistent collocation points}

First, suppose that
the assumptions in Theorem~\ref{Thm:Sol-Exist-Unique}
are fulfilled with $\domD=\SEt(\domD_d)$.
Then, approximating the equation $(I - \Fred)u = g$
by $(I - \FredSEn)\vSEn = g$
(i.e., applying Corollary~\ref{Cor:SE-Sinc-quad} as $\Fred\approx\FredSEn$),
set the approximate solution $\vSEn$ as
\begin{equation}
 \vSEn(t)
= g(t) + h \sum_{j=-N}^N k(t,\SEt(jh)) \vSEn(\SEt(jh))\SEtDiv(jh),
\label{Def:vSEn}
\end{equation}
where the mesh size $h$ is selected by the formula~\eqref{Def-h-SE}.
Setting $(2N+1)$ collocation points $t=\tSE_i$ $(i=-N,\,\ldots,\,N)$ as
\begin{equation}
 \tSE_i = \SEt(ih) \quad (i = -N,\,\ldots,\,N),
\label{eq:consistent-SE-sample}
\end{equation}
we obtain a $(2N+1)\times(2N+1)$ system of linear equations:
\begin{align}
 \sum_{j=-N}^N
  \left\{\delta_{ij}- h k(\tSE_i,\tSE_j)\SEtDiv(jh)\right\}\vSEn(\tSE_j)
&= g(\tSE_i),\quad i=-N,\ldots,\,N.
\label{Def:SE-Sinc-linear-eq-Nyst}
\end{align}
Let us define $m$ by $m=2N+1$, and let $I_m$ be
an $m\times m$ identity matrix,
and let $\tilde{K}_m^{\textSE}$ be an $m\times m$ matrix
whose $(i, j)$-th element is $h k(\tSE_{i},\tSE_{j})\SEtDiv(jh)$
$(i, j = -N,\,\ldots,\,N)$;
more precisely,
\[
 \tilde{K}_m^{\textSE}
=\begin{bmatrix}
  h k(\tSE_{-N}, \tSE_{-N})\SEtDiv(-Nh) & \cdots
& h k(\tSE_{-N}, \tSE_{j})\SEtDiv(jh) & \cdots
& h k(\tSE_{-N}, \tSE_{N})\SEtDiv(Nh) \\
\vdots& & \vdots & & \vdots \\
  h k(\tSE_{i}, \tSE_{-N})\SEtDiv(-Nh) & \cdots
& h k(\tSE_{i}, \tSE_{j})\SEtDiv(jh) & \cdots
& h k(\tSE_{i}, \tSE_{N})\SEtDiv(Nh) \\
\vdots & & \vdots & & \vdots \\
  h k(\tSE_{N}, \tSE_{-N})\SEtDiv(-Nh) & \cdots
& h k(\tSE_{N}, \tSE_{j})\SEtDiv(jh) & \cdots
& h k(\tSE_{N}, \tSE_{N})\SEtDiv(Nh) \\
 \end{bmatrix}.
\]
Furthermore, let
$\mathbd{\tilde{g}}_m^{\textSE}=[g(\tSE_{-N}),\,\ldots,\,g(\tSE_N)]^{\trans}$
and $\tilde{u}_j = \vSEn(jh)$ $(j=-N,\,\ldots,\,N)$.
Then,
the resulting system of linear equations~\eqref{Def:SE-Sinc-linear-eq-Nyst}
can be written in the matrix-vector form:
\begin{equation}
(I_m - \tilde{K}_m^{\textSE})\mathbd{\tilde{u}}_m
 = \mathbd{\tilde{g}}_m^{\textSE},
\label{Def:SE-Sinc-linear-eq-matvec-Nyst}
\end{equation}
where
$\mathbd{\tilde{u}}_m=[\tilde{u}_{-N},\,\ldots,\,\tilde{u}_N]^{\trans}$.
By the use of the solution vector $\mathbd{\tilde{u}}_m$,
let us introduce another approximate solution $\tilde{u}_N^{\textSE}$ as
\begin{equation}
 \tilde{u}_N^{\textSE}(t)
=\tilde{u}_{-N}\omega_a(t)
+\sum_{j=-N}^{N}\left(\tilde{u}_j
 - \tilde{u}_{-N}\frac{b-\SEt(jh)}{b-a}
 - \tilde{u}_{N} \frac{\SEt(jh)-a}{b-a}
\right)S(j,h)(\SEtInv(t))
+\tilde{u}_{N}\omega_b(t),
\label{Def:tilde-uSEn}
\end{equation}
where the mesh size $h$ is selected by the formula~\eqref{Def-h-SE}.
Between $\vSEn(t)$ and $\tilde{u}_N^{\textSE}(t)$,
the relationship $\tilde{u}_N^{\textSE} = \ProjSE\vSEn$ holds.
In summary, by solving the equation~\eqref{Def:SE-Sinc-linear-eq-matvec-Nyst}
for the coefficients $\mathbd{\tilde{u}}_m$,
the approximate solution $\tilde{u}_N^{\textSE}$
is determined by~\eqref{Def:tilde-uSEn}.
This is the new SE-Sinc-collocation method.
Its convergence theorem is written as follows,
which is proved in Section~\ref{sec:proof-SE}.

\begin{theorem}
\label{Thm:converge_new-SE}
Suppose that
the assumptions in Theorem~\ref{Thm:Sol-Exist-Unique} are fulfilled
with $\domD=\SEt(\domD_d)$ for $d\in(0,\pi)$.
Then, there exists a positive integer $N_0$ such that
the equation~\eqref{Def:SE-Sinc-linear-eq-matvec-Nyst}
is uniquely solvable for all $N\geq N_0$.
Furthermore, there exists a constant $C$ independent of $N$
such that for all $N\geq N_0$,
\[
\|u-\tilde{u}_N^{\textSE}\|_{C([a,b])}
 \leq C \sqrt{N}\exp\left(-\sqrt{\pi d \alpha N}\right).
\]
\end{theorem}

\subsection{New DE-Sinc-collocation method with consistent collocation points}

Here,
the new DE-Sinc-collocation method is described;
it replaces the SE transformation in the new SE-Sinc-collocation method
with the DE transformation.
Suppose that
the assumptions in Theorem~\ref{Thm:Sol-Exist-Unique}
are fulfilled with $\domD=\DEt(\domD_d)$.
Then, approximating the equation $(I - \Fred)u = g$
by $(I - \FredDEn)\vDEn = g$
(i.e., applying Corollary~\ref{Cor:DE-Sinc-quad} as $\Fred\approx\FredDEn$),
set the approximate solution $\vDEn$ as
\begin{equation}
 \vDEn(t)
= g(t) + h \sum_{j=-N}^N k(t,\DEt(jh)) \vDEn(\DEt(jh))\DEtDiv(jh),
\label{Def:vDEn}
\end{equation}
where the mesh size $h$ is selected by the formula~\eqref{Def-h-DE}.
Setting $(2N+1)$ collocation points $t=\tDE_i$ $(i=-N,\,\ldots,\,N)$ as
\begin{equation}
 \tDE_i = \DEt(ih) \quad (i = -N,\,\ldots,\,N),
\label{eq:consistent-DE-sample}
\end{equation}
we obtain an $m\times m$ system of linear equations (recall $m=2N+1$):
\begin{equation}
(I_m - \tilde{K}_m^{\textDE})\mathbd{\tilde{u}}_m
 = \mathbd{\tilde{g}}_m^{\textDE},
\label{Def:DE-Sinc-linear-eq-matvec-Nyst}
\end{equation}
where
$\mathbd{\tilde{g}}_m^{\textDE}=[g(\tDE_{-N}),\,\ldots,\,g(\tDE_N)]^{\trans}$,
and $\tilde{K}_m^{\textDE}$ is an $m\times m$ matrix
whose $(i, j)$-th element is $h k(\tDE_{i},\tDE_{j})\DEtDiv(jh)$
$(i, j = -N,\,\ldots,\,N)$;
more precisely,
\[
 \tilde{K}_m^{\textDE}
=\begin{bmatrix}
  h k(\tDE_{-N}, \tDE_{-N})\DEtDiv(-Nh) & \cdots
& h k(\tDE_{-N}, \tDE_{j})\DEtDiv(jh) & \cdots
& h k(\tDE_{-N}, \tDE_{N})\DEtDiv(Nh) \\
\vdots& & \vdots & & \vdots \\
  h k(\tDE_{i}, \tDE_{-N})\DEtDiv(-Nh) & \cdots
& h k(\tDE_{i}, \tDE_{j})\DEtDiv(jh) & \cdots
& h k(\tDE_{i}, \tDE_{N})\DEtDiv(Nh) \\
\vdots & & \vdots & & \vdots \\
  h k(\tDE_{N}, \tDE_{-N})\DEtDiv(-Nh) & \cdots
& h k(\tDE_{N}, \tDE_{j})\DEtDiv(jh) & \cdots
& h k(\tDE_{N}, \tDE_{N})\DEtDiv(Nh) \\
 \end{bmatrix}.
\]
By the use of the solution vector $\mathbd{\tilde{u}}_m$,
let us introduce another approximate solution $\tilde{u}_N^{\textDE}$ as
\begin{equation}
 \tilde{u}_N^{\textDE}(t)
=\tilde{u}_{-N}\omega_a(t)
+\sum_{j=-N}^{N}\left(\tilde{u}_j
 - \tilde{u}_{-N}\frac{b-\DEt(jh)}{b-a}
 - \tilde{u}_{N} \frac{\DEt(jh)-a}{b-a}
\right)S(j,h)(\DEtInv(t))
+\tilde{u}_{N}\omega_b(t),
\label{Def:tilde-uDEn}
\end{equation}
where the mesh size $h$ is selected by the formula~\eqref{Def-h-DE}.
Between $\vDEn(t)$ and $\tilde{u}_N^{\textDE}(t)$,
the relationship $\tilde{u}_N^{\textDE} = \ProjDE\vDEn$ holds.
In summary, by solving the equation~\eqref{Def:DE-Sinc-linear-eq-matvec-Nyst}
for the coefficients $\mathbd{\tilde{u}}_m$,
the approximate solution $\tilde{u}_N^{\textDE}$
is determined by~\eqref{Def:tilde-uDEn}.
This is the new DE-Sinc-collocation method.
Its convergence theorem is written as follows,
which is proved in Section~\ref{sec:proof-DE}.

\begin{theorem}
\label{Thm:converge_new-DE}
Suppose that
the assumptions in Theorem~\ref{Thm:Sol-Exist-Unique} are fulfilled
with $\domD=\DEt(\domD_d)$ for $d\in(0,\pi/2)$.
Then, there exists a positive integer $N_0$ such that
the equation~\eqref{Def:DE-Sinc-linear-eq-matvec-Nyst}
is uniquely solvable for all $N\geq N_0$.
Furthermore, there exists a constant $C$ independent of $N$
such that for all $N\geq N_0$,
\[
\|u-\tilde{u}_N^{\textDE}\|_{C([a,b])}
\leq C \exp\left\{\frac{-\pi d N}{\log(2 \pi d /\alpha)}\right\}.
\]
\end{theorem}

\subsection{Discussion on the differences between the original and new Sinc-collocation methods}
\label{sec:discussion}

The primary difference between the original and new Sinc-collocation methods
lies in the system of linear equations to be solved.
Recall that
the system of equations in the original methods is~\eqref{Def:SE-Sinc-linear-eq-matvec}
or~\eqref{Def:DE-Sinc-linear-eq-matvec},
whereas
that in the new methods is~\eqref{Def:SE-Sinc-linear-eq-matvec-Nyst}
or~\eqref{Def:DE-Sinc-linear-eq-matvec-Nyst}.
In the case of the original methods,
users must handle the exceptional cases at
$i,\,j = \pm (N+1)$,
due to the inconsistency of the collocation
points~\eqref{Def:SE-Sinc-Sampling}
or~\eqref{Def:DE-Sinc-Sampling}.
In the case of the new methods,
the systems of linear equations are considerably simple and easy to implement,
owing to the consistent collocation
points~\eqref{eq:consistent-SE-sample}
or~\eqref{eq:consistent-DE-sample}.
Especially in the matrices $\kSEn$ and $\kDEn$ of the original methods,
the elements in $j=\pm (N+1)$ include a summation term $\FredSEn$ or $\FredDEn$,
which may incur a high computational cost.
In contrast, such a term is not included in the matrices
$\tilde{K}_m^{\textSE}$ and $\tilde{K}_m^{\textDE}$ of the new methods.
For this reason,
the difference in the computation time between the original and new methods
is significant in the case where
$k$ is computationally complex for evaluation.

One minor difference is that
the size of the systems~\eqref{Def:SE-Sinc-linear-eq-matvec-Nyst}
and~\eqref{Def:DE-Sinc-linear-eq-matvec-Nyst}
is $m=2N+1$,
which is slightly smaller than
the size of the systems~\eqref{Def:SE-Sinc-linear-eq-matvec}
and~\eqref{Def:DE-Sinc-linear-eq-matvec}: $n=2N+3$.
It also should be noted that the original and new methods
share exactly the same basis functions.
Accordingly,
the number of the basis functions is the same: $n=2N+3$.
In terms of the convergence rates,
as Theorems~\ref{Thm:converge_RandZ}--\ref{Thm:converge_DE-Sinc}
and Theorems~\ref{Thm:converge_new-SE}--\ref{Thm:converge_new-DE} show,
there is no difference
between the original and new methods.

In short, the new methods have an advantage in terms of practicality of implementation
and computational cost,
without the loss of convergence performance.
This is also confirmed in the next section.

\section{Numerical examples}
\label{sec:numer}

In this section, numerical comparisons
between the four methods
(original SE/DE-Sinc-collocation methods
and new SE/DE-Sinc-collocation methods)
are presented.
The computation was performed on MacBook Air
with 1.7 GHz Intel Core i7 with 8GB memory, running
Mac OS X 10.12.6.
The computation programs were implemented in the C programming language
with double-precision
floating-point arithmetic, and compiled with Apple LLVM version 9.0.0
with the ``\texttt{-O2}'' option.
LAPACK in Apple's Accelerate framework was used
for computation of the system of linear equations.
The source code for all programs is available at
\url{https://github.com/okayamat/sinc-colloc-fredholm}.

We consider the following four examples.

\begin{example}[Delves and Mohamed~{\cite[special case of Family 4 in Chap.~10]{delves85}}]
\label{ex:1}
Consider the following equation
\[
 u(t) - \int_0^{1} t s u(s)\dd s
= \frac{r}{(t - 1/2)^2 + r^2} - t\arctan\left(\frac{1}{2r}\right),
\quad 0\leq t\leq 1,
\]
whose solution is $u(t)=r/\{(t - 1/2)^2 + r^2\}$.
For the experiment, $r$ was set as $r=1/2$.
\end{example}

\begin{example}[Delves and Mohamed~{\cite[Example 4.2.5b]{delves85}}]
\label{ex:2}
Consider the following equation
\[
 u(t) - \int_0^{\pi/2} (t s)^{3/4} u(s)\dd s
= t^{1/2}\left\{
1 - \frac{\pi^2}{9}\left(\frac{\pi t}{2}\right)^{1/4}
\right\},\quad 0\leq t\leq \pi/2,
\]
whose solution is $u(t)=t^{1/2}$.
\end{example}

\begin{example}[Okayama et al.~{\cite[Example 9.4]{OMS}}]
\label{ex:3}
Consider the following equation
\[
 u(t) - \int_0^{1} t^{\sqrt{3} - 1}
\sum_{l=1}^{100} s^{a_l}(1 - s)^{1 - b_l} u(s)\dd s
= t^{1/2} - t^{\sqrt{3} - 1}
\sum_{l=1}^{100} B(a_l + 3/2, 2 - b_l),\quad 0\leq t\leq 1,
\]
whose solution is $u(t)=t^{1/2}$.
Here, $\Beta(p,q)$ is the beta function,
$a_l = (3/\pi)^l$, and $b_l = (2\sqrt{2}/3)^{l}$.
\end{example}

\begin{example}[Okayama et al.~{\cite[Example 9.5]{OMS}}]
\label{ex:4}
Consider the following equation
\begin{align*}
& u(t) - 2\int_{-1}^{1}
\frac{(1 - s^2)^{(2-t^2)/(2+t^2)}}{(2+t^2)(1+s^{20})}
\left\{5(2+t^2)s^{18}(1-s^2)+(s^{20}+1)(s^{21}+s+2)\right\}
u(s)\dd s\\
&= \frac{2t}{1+t^{20}} - \frac{4}{2+t^2}\Beta(3/2,4/(2+t^2)),
\quad -1\leq t\leq 1,
\end{align*}
whose solution is $u(t)=2t/(1+t^{20})$.
\end{example}

Following Okayama et al.~\cite{OMS},
the values of the parameters $\alpha$ and $d$ were chosen
as shown in Table~\ref{tbl:param}.
Here, $\pim=3.14$ was used as a slightly smaller value than $\pi$.
Because there is no difference between the assumptions
in the theorems of the original and new methods,
the methods share the same values for those parameters.
The number of lines of the main programs
is shown in Table~\ref{tbl:program}.
The table
reveals the advantage of the new SE/DE-Sinc collocation methods
in terms of implementation.
Furthermore, note that the differences
between the four examples are small,
even though Examples~\ref{ex:3} and~\ref{ex:4} are
extremely difficult for Gaussian quadrature-based methods.

\begin{table}
\caption{Values of the parameters $\alpha$ and $d$. $\pim=3.14$.}
\label{tbl:param}
\centering
\begin{tabular}{ccccc} \hline
 & Example~\ref{ex:1} & Example~\ref{ex:2}
 & Example~\ref{ex:3} & Example~\ref{ex:4} \\ \hline
Original SE, $\alpha$ & $1.0$    & $0.5$    & $0.5$    & $1.0$ \\
New SE, $\alpha$ & $1.0$    & $0.5$    & $0.5$    & $1.0$ \\
Original DE, $\alpha$ & $1.0$    & $0.5$    & $0.5$    & $1.0$ \\
New DE, $\alpha$ & $1.0$    & $0.5$    & $0.5$    & $1.0$ \\
Original SE, $d$      & $\pim/2$ & $\pim$   & $\pim$   & $\pim/2$ \\
New SE, $d$      & $\pim/2$ & $\pim$   & $\pim$   & $\pim/2$ \\
Original DE, $d$      & $\pim/6$ & $\pim/2$ & $\pim/2$ & $0.125$ \\
New DE, $d$      & $\pim/6$ & $\pim/2$ & $\pim/2$ & $0.125$ \\
\hline
\end{tabular}
\end{table}
\begin{table}
\caption{Number of lines of the main programs.}
\label{tbl:program}
\centering
\begin{tabular}{ccccc} \hline
 & Example~\ref{ex:1} & Example~\ref{ex:2}
 & Example~\ref{ex:3} & Example~\ref{ex:4} \\ \hline
Original SE & $171$ & $169$ & $204$ & $182$ \\
Original DE & $171$ & $169$ & $204$ & $183$ \\
New SE      & $127$ & $125$ & $160$ & $138$ \\
New DE      & $127$ & $125$ & $160$ & $139$ \\
\hline
\end{tabular}
\end{table}

\begin{figure}
\begin{center}
 \begin{minipage}{0.45\linewidth}
  \includegraphics[width=\linewidth]{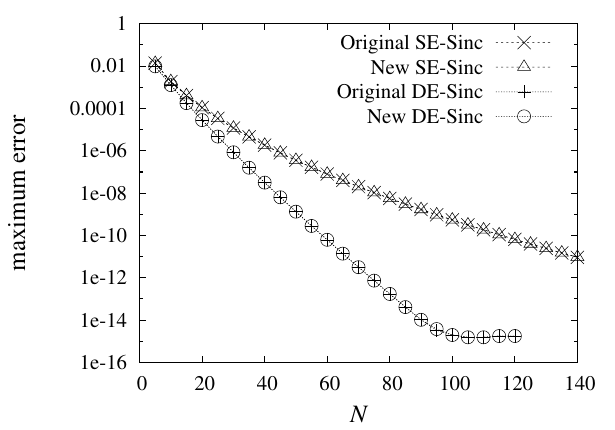}
  \caption{Convergence with respect to $N$ for Example~\ref{ex:1}.}
  \label{Fig:Ex1error}
 \end{minipage}
 \begin{minipage}{0.01\linewidth}
 \mbox{ }
 \end{minipage}
 \begin{minipage}{0.45\linewidth}
  \includegraphics[width=\linewidth]{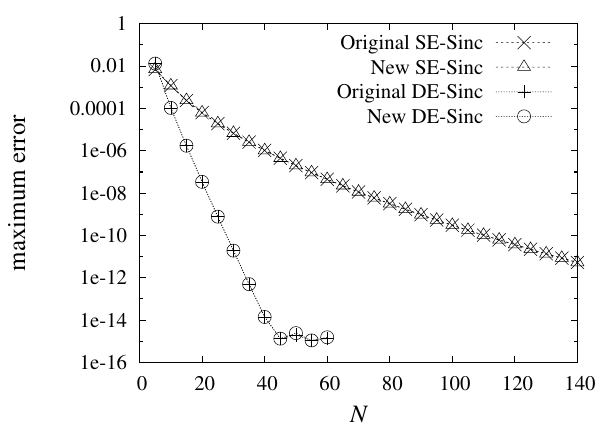}
  \caption{Convergence with respect to $N$ for Example~\ref{ex:2}.}
  \label{Fig:Ex2error}
 \end{minipage}
\end{center}
\begin{center}
 \begin{minipage}{0.45\linewidth}
  \includegraphics[width=\linewidth]{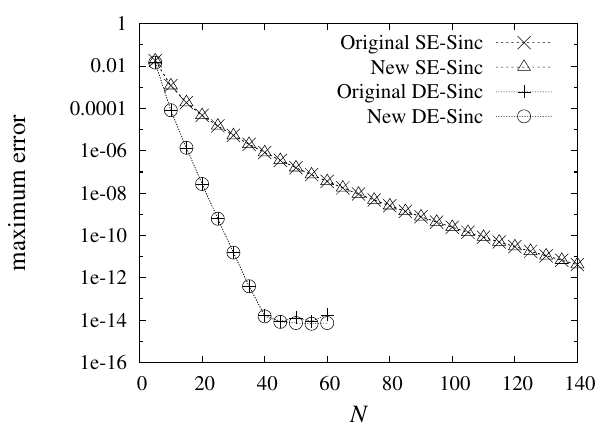}
  \caption{Convergence with respect to $N$ for Example~\ref{ex:3}.}
  \label{Fig:Ex3error}
 \end{minipage}
 \begin{minipage}{0.01\linewidth}
 \mbox{ }
 \end{minipage}
 \begin{minipage}{0.45\linewidth}
  \includegraphics[width=\linewidth]{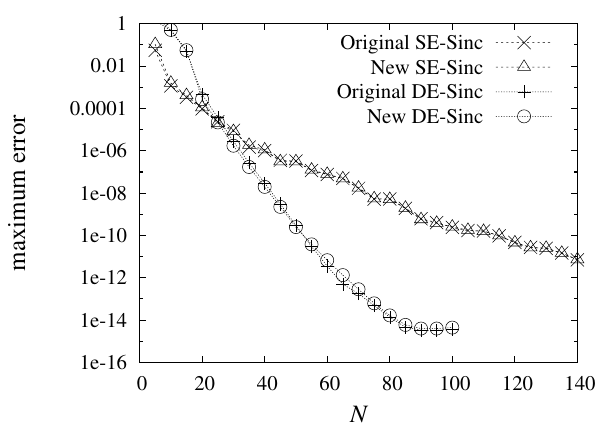}
  \caption{Convergence with respect to $N$ for Example~\ref{ex:4}.}
  \label{Fig:Ex4error}
 \end{minipage}
\end{center}
\end{figure}

\begin{figure}
\begin{center}
 \begin{minipage}{0.45\linewidth}
  \includegraphics[width=\linewidth]{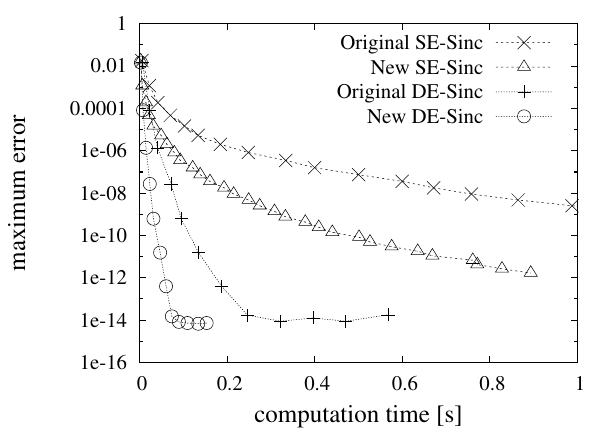}
  \caption{Convergence with respect to computation time for Example~\ref{ex:3}.}
  \label{Fig:Ex3time}
 \end{minipage}
\end{center}
\end{figure}

Figures~\ref{Fig:Ex1error}--\ref{Fig:Ex4error} show
the convergence profile with respect to $N$.
We observe that the original and new SE-Sinc-collocation methods
converge at the same rate, $\OO(\exp(-c\sqrt{N}))$,
which is suggested by
Theorems~\ref{Thm:converge_RandZ} and~\ref{Thm:converge_new-SE}.
We also observe that the original and new DE-Sinc-collocation methods
converge at the same rate, $\OO(\exp(-c N/\log N))$,
which is suggested by
Theorems~\ref{Thm:converge_DE-Sinc} and~\ref{Thm:converge_new-DE}.
However, especially in Example~\ref{ex:3},
the notable advantage of the new SE/DE-Sinc-collocation methods
may be observed in terms of the computational cost.
The convergence profile with respect to the computation time
is displayed in Figure~\ref{Fig:Ex3time}.
The difference in the performance
is due to the functions $k$ and $g$ in Example~\ref{ex:3},
which are relatively computationally complex for evaluation.
In other examples, such a significant difference was not observed.


\section{Proof of convergence theorems}
\label{sec:proofs}

This section provides proofs for
Theorems~\ref{Thm:converge_new-SE} and~\ref{Thm:converge_new-DE}.
In this section, we write $X = C([a,b])$ for short.

\subsection{Proof of Theorem~\ref{Thm:converge_new-SE}}
\label{sec:proof-SE}

Recall that by solving~\eqref{Def:SE-Sinc-linear-eq-matvec-Nyst},
$\vSEn$ is determined by~\eqref{Def:vSEn},
and $\tilde{u}_N^{\textSE}$ is determined by~\eqref{Def:tilde-uSEn}.
The key point of the proof of Theorem~\ref{Thm:converge_new-SE}
is the relation between those two approximate
solutions: $\tilde{u}_N^{\textSE}=\ProjSE\vSEn$.
Using the relation, we have
\begin{equation}
 \|u - \tilde{u}_N^{\textSE}\|_{X}
=\|u - \ProjSE \vSEn\|_{X}
\leq
 \|u - \ProjSE u\|_{X}
+\|\ProjSE u - \ProjSE \vSEn\|_{X}
\leq
 \|u - \ProjSE u\|_{X}
+\|\ProjSE\|_{\mathcal{L}(X,X)} \|u - \vSEn\|_{X}.
\label{eq:SE-first-estimate}
\end{equation}
Under the assumptions of Theorem~\ref{Thm:converge_new-SE},
$u\in\MC_{\alpha}(\SEt(\domD_d))$ holds by Theorem~\ref{Thm:Sol-Exist-Unique}.
Thus, using Theorem~\ref{Thm:general-SE-Sinc},
we can evaluate the first term in~\eqref{eq:SE-first-estimate} as
\[
 \|u - \ProjSE u\|_{X} \leq C_1 \sqrt{N} \exp\left(-\sqrt{\pi d \alpha N}\right)
\]
for a constant $C_1$ independent of $N$.
As for the second term in~\eqref{eq:SE-first-estimate},
the following result was obtained.

\begin{theorem}[Okayama et al.~{\cite[Theorem~6.3]{OMS}}]
 \label{Thm:estimate_u-v}
Suppose that
the assumptions in Theorem~\ref{Thm:Sol-Exist-Unique} are fulfilled
with $\domD=\SEt(\domD_d)$ for $d\in(0,\,\pi)$.
Then, there exists a positive integer $N_0$ such that
for all $N\geq N_0$,
the equation~\eqref{Def:vSEn} has a unique solution $\vSEn\in X$.
Furthermore, there exists a constant $C_2$ independent of $N$
such that for all $N\geq N_0$,
\[
 \|u-\vSEn\|_X \leq C_2 \exp\left(-\sqrt{\pi d \alpha N}\right).
\]
\end{theorem}

In addition to providing the estimate of $\|u-\vSEn\|_X$,
this theorem guarantees the existence of the unique solution $\vSEn$;
this means that the system~\eqref{Def:SE-Sinc-linear-eq-matvec-Nyst}
is uniquely solvable.
Furthermore, we can estimate the remaining term
$\|\ProjSE\|_{\mathcal{L}(X,X)}$ as follows.

\begin{lemma}
\label{lem:ProjSE-estimate}
There exists a constant $C_3$ independent of $N$
such that $\|\ProjSE\|_{\mathcal{L}(X,X)}\leq C_3\log(N+1)$.
\end{lemma}

The lemma is deduced immediately by using the following result.

\begin{lemma}[Stenger~{\cite[p.\,142]{Ste1}}]
\label{Lem:Sinc-Real-Sum}
Let $h>0$. Then, it holds that
\[
\sup_{x\in\mathbb{R}}\sum_{j=-N}^N |S(j,h)(x)|
\leq \frac{2}{\pi}(3 + \log N).
\]
\end{lemma}

Thus,
we can evaluate the second term in~\eqref{eq:SE-first-estimate} as
\[
\|\ProjSE\|_{\mathcal{L}(X,X)}\|u - \ProjSE\vSEn\|_{X}
\leq C_3 \log(N+1) \cdot C_2 \exp\left(-\sqrt{\pi d \alpha N}\right).
\]
Summing up the above results, we establish Theorem~\ref{Thm:converge_new-SE}.

\subsection{Proof of Theorem~\ref{Thm:converge_new-DE}}
\label{sec:proof-DE}

Recall that by solving~\eqref{Def:DE-Sinc-linear-eq-matvec-Nyst},
$\vDEn$ is determined by~\eqref{Def:vDEn},
and $\tilde{u}_N^{\textDE}$ is determined by~\eqref{Def:tilde-uDEn}.
In the proof of Theorem~\ref{Thm:converge_new-DE} as well,
the key point is the relation between those two approximate
solutions: $\tilde{u}_N^{\textDE}=\ProjDE\vDEn$.
Using the relation, we have
\begin{equation}
 \|u - \tilde{u}_N^{\textDE}\|_{X}
=\|u - \ProjDE \vDEn\|_{X}
\leq
 \|u - \ProjDE u\|_{X}
+\|\ProjDE u - \ProjDE \vDEn\|_{X}
\leq
 \|u - \ProjDE u\|_{X}
+\|\ProjDE\|_{\mathcal{L}(X,X)} \|u - \vDEn\|_{X}.
\label{eq:DE-first-estimate}
\end{equation}
Under the assumptions of Theorem~\ref{Thm:converge_new-DE},
$u\in\MC_{\alpha}(\DEt(\domD_d))$ holds by Theorem~\ref{Thm:Sol-Exist-Unique}.
Thus, using Theorem~\ref{Thm:general-DE-Sinc},
we can evaluate the first term in~\eqref{eq:DE-first-estimate} as
\[
 \|u - \ProjDE u\|_{X} \leq C_1 \exp\left\{\frac{-\pi d N}{\log(2 d N/\alpha)}\right\}
\]
for a constant $C_1$ independent of $N$.
As for the second term in~\eqref{eq:DE-first-estimate},
the following result was obtained.

\begin{theorem}[Okayama et al.~{\cite[Theorem~8.1]{OMS}}]
 \label{Thm:estimate_u-v-DE}
Suppose that
the assumptions in Theorem~\ref{Thm:Sol-Exist-Unique} are fulfilled
with $\domD=\DEt(\domD_d)$ for $d\in(0,\,\pi/2)$.
Then, there exists a positive integer $N_0$ such that
for all $N\geq N_0$,
the equation~\eqref{Def:vDEn} has a unique solution $\vDEn\in X$.
Furthermore, there exists a constant $C_2$ independent of $N$
such that for all $N\geq N_0$,
\[
 \|u-\vDEn\|_X \leq C_2 \exp\left\{\frac{-2\pi d N}{\log(2 d N/\alpha)}\right\}.
\]
\end{theorem}

In addition to providing the estimate of $\|u-\vDEn\|_X$,
this theorem guarantees the existence of the unique solution $\vDEn$;
this means that the system~\eqref{Def:DE-Sinc-linear-eq-matvec-Nyst}
is uniquely solvable.
Furthermore,
in the same way as Lemma~\ref{lem:ProjSE-estimate},
we can estimate the remaining term
$\|\ProjDE\|_{\mathcal{L}(X,X)}$ as follows.

\begin{lemma}
There exists a constant $C_3$ independent of $N$
such that $\|\ProjDE\|_{\mathcal{L}(X,X)}\leq C_3\log(N+1)$.
\end{lemma}

Thus,
we can evaluate the second term in~\eqref{eq:DE-first-estimate} as
\[
\|\ProjDE\|_{\mathcal{L}(X,X)}\|u - \ProjDE\vDEn\|_{X}
\leq C_3 \log(N+1) \cdot C_2 \exp\left\{\frac{-2\pi d N}{\log(2 d N/\alpha)}\right\}.
\]
Summing up the above results, we establish Theorem~\ref{Thm:converge_new-DE}.

\section{Concluding remarks}
\label{sec:conclusion}

Okayama et al.~\cite{OMS} presented
the SE/DE-Sinc-collocation methods for the equation~\eqref{eq:Fredholm}.
Those methods can attain exponential convergence even if
the functions $k$ and $g$ have derivative singularity at the endpoints.
However, those methods are challenging or painful to implement
owing to the inconsistent collocation points~\eqref{Def:SE-Sinc-Sampling}
or~\eqref{Def:DE-Sinc-Sampling}.
To address the situation, this paper proposed
new SE/DE-Sinc-collocation methods
that employ consistent collocation points~\eqref{eq:consistent-SE-sample}
or~\eqref{eq:consistent-DE-sample}.
The results of a theoretical analysis proved that the new SE/DE-Sinc-collocation methods
exhibit the same convergence property
as the original SE/DE-Sinc-collocation methods.
As shown in Section~\ref{sec:numer},
compared with the original methods,
the number of lines of the programs in the new methods was reduced to 70--80~\%.
Furthermore, the advantage of the new methods in terms of the computational cost
may be observed in the equations with computationally complex functions $k$ and $g$.
It should also be emphasized that
the condition numbers of the resulting matrices
in the new methods
($I_m - \tilde{K}_m^{\textSE}$ and $I_m - \tilde{K}_m^{\textDE}$)
are good~\cite{OMS3},
although such a result is not obtained yet
for the original methods.

The followings issues may be considered in future studies.
First, a Sinc-collocation method for Volterra integral equations
was proposed by Rashidinia and Zarebnia~\cite{RZ3},
and results similar to those of the present work can be established for that method.
Second, the error analysis for the Sinc-Nystr\"{o}m methods
can be improved compared with Theorems~\ref{Thm:estimate_u-v}
and~\ref{Thm:estimate_u-v-DE}.
As described in Remark~\ref{rem:choice-h},
the mesh size $h$ for the SE/DE-Sinc quadrature
is not selected optimally because it can share the same $h$
as in the SE/DE-Sinc approximation, which is also used in
the Sinc-collocation methods.
However, in the case of the Sinc-Nystr\"{o}m methods,
only the SE/DE-Sinc quadrature is used,
and the mesh size $h$ can be optimally chosen.
The error of this case should be rigorously analyzed for both cases of the SE transformation~\cite{RZ2}
and DE transformation~\cite{MM}.
Third, explicit error bounds are desired for verified numerical computation.
If we can compute the constant $C$ in Theorem~\ref{Thm:converge_new-SE}
or~\ref{Thm:converge_new-DE},
we can guarantee the accuracy of the method by computing
the right-hand side of the inequality.
This can be performed by examining the constants appearing in the
convergence theorems, and using the results for basic
approximation formulas~\cite{OMS4}.
We plan to report further results on these issues in future work.



\begin{thebibliography}{99}
 \bibitem{delves85} L. M. Delves and J. L. Mohamed.
\newblock {\em Computational {M}ethods for {I}ntegral {E}quations}.
\newblock Cambridge University Press, Cambridge, 1985.
 \bibitem{MS} M. Mori and M. Sugihara.
\newblock The double-exponential transformation in numerical analysis.
\newblock {\em J. Comput.\ Appl.\ Math.}, 127(1--2):287--296, 2001.
 \bibitem{MM} M. Muhammad, A. Nurmuhammad, M. Mori and M. Sugihara.
\newblock Numerical solution of integral equations by means of the {S}inc
 collocation method based on the double exponential transformation.
\newblock {\em J. Comput.\ Appl.\ Math.}, 177(2):269--286, 2005.
 \bibitem{TO2} T. Okayama.
\newblock A note on the {S}inc approximation with boundary treatment.
\newblock {\em JSIAM Lett.}, 5:1--4, 2013.
 \bibitem{TO} T. Okayama.
\newblock Theoretical analysis of {S}inc-collocation methods and
 {S}inc-{N}ystr\"{o}m  methods for systems of initial value problems.
\newblock {\em BIT Numer.\ Math.}, 58(1):199--220, 2018.
 \bibitem{OMS4}  T. Okayama, T. Matsuo and M. Sugihara.
\newblock Error estimates with explicit constants for {S}inc approximation,
 {S}inc quadrature and {S}inc indefinite integration.
\newblock {\em Numer.\ Math.}, 124(2):361--394, 2013.
 \bibitem{OMS}  T. Okayama, T. Matsuo and M. Sugihara.
\newblock Improvement of a {S}inc-collocation method for {F}redholm integral
 equations of the second kind.
\newblock {\em BIT Numer.\ Math.}, 51(2):339--366, 2011.
 \bibitem{OMS3}  T. Okayama, T. Matsuo and M. Sugihara.
\newblock On boundedness of the condition number of the coefficient matrices
 appearing in {S}inc-{N}ystr\"{o}m methods for {F}redholm integral equations
 of the second kind.
\newblock {\em JSIAM Letters}, 3:81--84, 2011.
 \bibitem{OMS2}  T. Okayama, T. Matsuo and M. Sugihara.
\newblock Sinc-collocation methods for weakly singular {F}redholm integral
 equations of the second kind.
\newblock {\em J. Comput.\ Appl.\ Math.}, 234(4):1211--1227, 2010.
 \bibitem{RZ2} J. Rashidinia and M. Zarebnia.
\newblock Convergence of approximate solution of system of {F}redholm
 integral equations.
\newblock {\em J. Math.\ Anal.\ Appl.}, 333(2):1216--1227, 2007.
 \bibitem{RZ1} J. Rashidinia and M. Zarebnia.
\newblock Numerical solution of linear integral equations by using
 {S}inc-collocation method.
\newblock {\em Appl.\ Math.\ Comput.}, 168(2):806--822, 2005.
 \bibitem{RZ3} J. Rashidinia and M. Zarebnia.
\newblock Solution of a {V}olterra integral equation by the
 {S}inc-collocation method.
\newblock {\em J. Comput.\ Appl.\ Math.}, 206(2):801--813, 2007.
 \bibitem{Ste1} F. Stenger.
\newblock {\em Numerical {M}ethods {B}ased on {S}inc and {A}nalytic {F}unctions}.
\newblock Springer-Verlag, New York, 1993.
 \bibitem{Ste3} F. Stenger.
\newblock {\em Handbook of {S}inc {N}umerical {M}ethods}.
\newblock CRC Press, Boca Raton, FL, 2011.
 \bibitem{Ste2} F. Stenger.
\newblock Summary of {S}inc numerical methods.
\newblock {\em J. Comput.\ Appl.\ Math.}, 121(1--2):379--420, 2000.
 \bibitem{SugiMatsu} M. Sugihara and T. Matsuo.
\newblock Recent developments of the {S}inc numerical methods.
\newblock {\em J. Comput.\ Appl.\ Math.}, 164--165(1):673--689, 2004.
 \bibitem{tanaka09:_class_of} K. Tanaka, M. Sugihara and K. Murota.
\newblock Function classes for successful {DE}-{S}inc approximations.
\newblock {\em Math. Comput.}, 78(267):1553--1571, 2009.
\end{thebibliography}
\end{document}